\documentclass[11pt]{article}%
\usepackage{amsmath}
\usepackage{amsfonts}
\usepackage{amssymb}
\usepackage{graphicx}%
\begin{document}

\title{The Bramble-Hilbert Lemma\thanks{I place this document in the public domain
and you can do with it anything you want.}}
\author{Jan Mandel\\University of Colorado Denver}
\date{\relax}
\maketitle

\begin{abstract}
This is an introductory document surveying several results in polynomial
approximation, known as the Bramble-Hilbert lemma.

\end{abstract}

\section{Introduction}

In numerical analysis, the Bramble-Hilbert lemma bounds the error of an
approximation of a function $u$ by a polynomial of order at most $m-1$ in
terms of derivatives of $u$ of order $m$. Both the error of the approximation
and the derivatives of $u$ are measured by $L^{p}$ norms on a bounded domain
in $\mathbb{R}^{n}$. This is similar to classical numerical analysis, where,
for example, the error of interpolation $u$ on an interval by a linear
function (that is, approximation by a polynomial of order one) can be bounded
using the second derivative of $u$. The difference is that the Bramble-Hilbert
lemma applies in any number of dimensions, not just one dimension, and the
approximation error and the derivatives of $u$ are measured by more general
norms involving averages, not just the maximum norm.

Additional assumptions on the domain are needed for the Bramble-Hilbert lemma
to hold. Essentially, the boundary of the domain must be \textquotedblleft
reasonable\textquotedblright. For example, domains that have a spike or a slit
with zero angle at the tip are excluded. Domains that are reasonable enough
include all convex domains and Lipschitz domains, which includes all domains
with a continuously differentiable boundary.

The main use of the Bramble-Hilbert lemma is to prove bounds on the error of
interpolation of function $u$ by an operator that preserves polynomials of
order up to $m-1$, in terms of the derivatives of $u$ of order $m$. This is an
essential step in error estimates for the finite element method. The
Bramble-Hilbert lemma is applied there on the domain consisting of one element.

\section{The one dimensional case}

Before stating the lemma in full generality, it is useful to look at some
simple special cases. In one dimension and for a function $u$ that has $m$
derivatives on interval $\left(  a,b\right)  $, the lemma reduces to%
\[
\inf_{v\in P_{m-1}}\bigl\Vert u^{\left(  k\right)  }-v^{\left(  k\right)
}\bigr\Vert_{L^{p}\left(  a,b\right)  }\leq C\left(  m\right)  \left(
b-a\right)  ^{m-k}\bigl\Vert u^{\left(  m\right)  }\bigr\Vert_{L^{p}\left(
a,b\right)  },
\]
where $P_{m-1}$ is the space of all polynomials of order at most $m-1$.

In the case when $p=\infty$, $m=2$, $k=1$, and $u$ is twice differentiable,
this means that there exists a polynomial $v$ of degree one such that for all
$x\in\left(  a,b\right)  $,%
\[
\left\vert u\left(  x\right)  -v\left(  x\right)  \right\vert \leq C\left(
b-a\right)  ^{2}\sup_{\left(  a,b\right)  }\left\vert u^{\prime\prime
}\right\vert ,\text{\quad for all }x\in\left(  a,b\right)  .
\]
This inequality follows from the well-known error estimate for linear
interpolation by choosing $v$ as the linear interpolant of $u$.

\section{Statement of the lemma}

Suppose $\Omega$ is a bounded domain in $\mathbb{R}^{n}$, $n\geq1$, with
boundary $\partial\Omega$ and diameter $d$. $W_{p}^{k}(\Omega)$ is the Sobolev
space of all function $u$ on $\Omega$ with weak derivatives $D^{\alpha}u$ of
order $\left\vert \alpha\right\vert $ up to $k$ in $L^{p}(\Omega)$. Here,
$\alpha=\left(  \alpha_{1},\alpha_{2},\ldots,\alpha_{n}\right)  $ is a
multiindex, $\left\vert \alpha\right\vert =$ $\alpha_{1}+\alpha_{2}%
+\cdots+\alpha_{n}$, and $D^{\alpha}$ denotes the derivative $\alpha_{1}$
times with respect to $x_{1}$, $\alpha_{2}$ times with respect to $\alpha_{2}%
$, and so on. The Sobolev seminorm on $W_{p}^{m}(\Omega)$ consists of the
$L^{p}$ norms of the highest order derivatives,%
\[
\left\vert u\right\vert _{W_{p}^{m}(\Omega)}=\left(  \sum_{\left\vert
\alpha\right\vert =m}\left\Vert D^{\alpha}u\right\Vert _{L^{p}(\Omega)}%
^{p}\right)  ^{1/p}\text{ if }1\leq p<\infty
\]
and%
\[
\left\vert u\right\vert _{W_{\infty}^{m}(\Omega)}=\max_{\left\vert
\alpha\right\vert =m}\left\Vert D^{\alpha}u\right\Vert _{L^{\infty}(\Omega)}%
\]

$P_{k}$ is the space of all polynomials of order up to $k$ on $\mathbb{R}^{n}%
$. Note that $D^{\alpha}v=0$ for all $v\in P_{m-1}$. and $\left\vert
\alpha\right\vert =m$, so $\left\vert u+v\right\vert _{W_{p}^{m}(\Omega)}$ has
the same value for any $v\in P_{k-1}$.

\textbf{Lemma (Bramble and Hilbert)} Under additional assumptions on the
domain $\Omega$, specified below, there exists a constant $C=C\left(
m,\Omega\right)  $ independent of $p$ and $u$ such that for any $u\in
W_{p}^{k}(\Omega)$ there exists a polynomial $v\in P_{m-1}$ such that for all
$k=0,\ldots,m$,%
\[
\left\vert u-v\right\vert _{W_{p}^{k}(\Omega)}\leq Cd^{m-k}\left\vert
u\right\vert _{W_{p}^{m}(\Omega)}.
\]

\section{The original result}

The lemma was proved by Bramble and Hilbert \cite{Bramble-1970-ELF} under the
assumption that $\Omega$ satisfies the \emph{strong cone property;} that is,
there exists a finite open covering $\left\{  O_{i}\right\}  $ of
$\partial\Omega$ and corresponding cones $\{C_{i}\}$ with vertices at the
origin such that $x+C_{i}$ is contained in $\Omega$ for any $x$ $\in\Omega\cap
O_{i}$.

The statement of the lemma here is a simple rewriting of the right-hand
inequality stated in Theorem 1 in \cite{Bramble-1970-ELF}. The actual
statement in \cite{Bramble-1970-ELF} is that the norm of the factorspace
$W_{p}^{m}(\Omega)/P_{m-1}$ is equivalent to the $W_{p}^{m}(\Omega)$ seminorm.
The $W_{p}^{m}(\Omega)$ is not the usual one but the terms are scaled with $d$
so that the right-hand inequality in the equivalence of the seminorms comes
out exactly as in the statement here.

In the original result, the choice of the polynomial is not specified, and the
dependence of the constant on the domain $\Omega$ is not given either.

\section{A constructive form}

An alternative result was given by Dupont and Scott \cite{Dupont-1980-PAF}
under the assumption that the domain $\Omega$ is \emph{star-shaped;} that is,
there exists a ball $B$ such that for any $x\in\Omega$, the closed convex hull
of $\left\{  x\right\}  \cup B$ is a subset of $\Omega$. Suppose that
$\rho_{\max}$ is the supremum of the diameters of such balls. The ratio
$\gamma=d/\rho_{\max}$ is called the chunkiness of $\Omega$.

Given a fixed ball $B$ as above, and a function $u$, the averaged Taylor
polynomial $Q^{m}u$ is defined as%
\[
Q^{m}u=\int\limits_{B}T_{y}^{m}u\left(  x\right)  \psi\left(  y\right)  dx,
\]
where
\[
T_{y}^{m}u\left(  x\right)  =\sum\limits_{k=0}^{m-1}\sum\limits_{\left\vert
\alpha\right\vert =k}\frac{1}{\alpha!}D^{\alpha}u\left(  y\right)  \left(
x-y\right)  ^{\alpha}%
\]
is the Taylor polynomial of degree at most $m-1$ of $u$ centered at $y$
evaluated at $x$, and $\psi\geq0$ is a function that has derivatives of all
orders, equals to zero outside of $B$, and such that%
\[
\int\limits_{B}\psi dx=1.
\]
Such function $\psi$ always exists.

Then the lemma holds with the constant $C=C\left(  m,n,\gamma\right)  $, that
is, the constant depends on the domain $\Omega$ only through its chunkiness
$\gamma$ and the dimension of the space $n$. For more details and a tutorial
treatment, see the monograph by Brenner and Scott \cite{Brenner-2002-MTF}. The
result can be extended to the case when the domain $\Omega$ is the union of a
finite number of star-shaped domains and more general polynomial spaces than
the space of all polynomials up to a given degree \cite{Dupont-1980-PAF}.

\section{Bound on linear functionals}

This result follows immediately from the above lemma, and it is also called
sometimes the Bramble-Hilbert lemma, for example by Ciarlet
\cite{Ciarlet-2002-FEM}. It is essentially Theorem 2 from
\cite{Bramble-1970-ELF}.

\textbf{Lemma} Suppose that $\ell$ is a continuous linear functional on
$W_{p}^{m}(\Omega)$ and $\left\Vert \ell\right\Vert _{W_{p}^{m}(\Omega
)^{^{\prime}}}$ its dual norm. Suppose that $\ell\left(  v\right)  =0$ for all
$v\in P_{m-1}$. Then there exists a constant $C=C\left(  \Omega\right)  $ such
that%
\[
\left\vert \ell\left(  u\right)  \right\vert \leq C\left\Vert \ell\right\Vert
_{W_{p}^{m}(\Omega)^{^{\prime}}}\left\vert u\right\vert _{W_{p}^{m}(\Omega)}.
\]

\bibliographystyle{plain}
\bibliography{../../bddc/bibliography/bddc}

\end{document}